# CAUCHY PROBLEM IN THE NON-CLASSICAL TREATMENT FOR ONE PSEUDOPARABOLIC EQUATION

**I.G. Mamedov**

*A.I. Huseynov Institute of Cybernetics of NAS of Azerbaijan. Az 1141,
Azerbaijan, Baku st. B. Vahabzade, 9
E-mail: ilgar-mammadov@rambler.ru*

## Abstract

*In the paper, we consider the Cauchy problem for a fifth order pseudoparabolic equation that appears in studying the issues of fluid filtration in fissured media, the moisture transfer in soils and etc. The Cauchy problem with non-classic conditions not requiring the agreement conditions are studied for a discontinuous coefficient pseudoparabolic equation. The equivalence of these conditions with the Cauchy classic condition is substantiated in the case when the solution of the stated problem is sought in S.L. Sobolev anisotropic space.*

## Problem statement

Let $G = G_1 \times G_2$, $G_k = (0, h_k)$, $k = \overline{1,2}$; $W_p^{(3,2)}(G)$, $1 \le p \le \infty$ be a space of all functions $u \in L_p(G)$ having the generalized derivatives, in S.L. Sobolev's sense $D_x^i D_y^j u \in L_p(G)$, $i = \overline{0,3}$; $j = \overline{0,2}$, where $D_t = \partial/\partial_t$ is a generalized differentiation operator in S.L. Sobolev's sense, $D_t^0$ is an identical transformation operator. We'll define the norm in S.L. Sobolev anisotropic space $W_p^{(3,2)}(G)$ by the following equality

$$\|u\|_{W_p^{(3,2)}(G)} = \sum_{i=0}^{3} \sum_{j=0}^{2} \left\| D_x^i D_y^j u \right\|_{L_p(G)}$$

Consider the equation

$$(V_{3,2}u)(x, y) \equiv D_x^3 D_y^2 u(x, y) + a_{2,2}(x, y)D_x^2 D_y^2 u(x, y) +$$
$$+ a_{3,1}(x, y)D_x^3 D_y u(x, y) + a_{2,1}(x, y)D_x^2 D_y u(x, y) +$$
$$+ a_{1,2}(x, y)D_x D_y^2 u(x, y) + a_{3,0}(x, y)D_x^3 u(x, y) + \quad (1)$$
$$+ \sum_{\substack{i=0 \\ i+j<3}}^{2} \sum_{j=0}^{2} a_{i,j}(x, y)D_x^i D_y^j u(x, y) = Z_{3,2}(x, y), \ (x, y) \in G,$$

where $u(x, y)$ is a desired function from $W_p^{(3,2)}(G)$.

Some classes of boundary value problems for equation (1) in definite sense are stated similar to the known boundary value problems for the parabolic equation $D_y u(x, y) = D_x^2 u(x, y)$. Therefore, many authors call the equation of the form (1) a pseudoparabolic equation. Notice that the equation under consideration is a generalization of many model equations of some processes (for example, of heat conductivity equation, string vibration equation, telegraph equation, generalized moisture transfer equation, Boussenesq-Liav equation and etc.).

In particular, many processes arising in theory of fluid filtration in fissured media [1-2] are described by pseudopa-rabolic equations with discontinuous coefficients [3-4] Such equations arise while describing a lot of real processes occurring in nature and engineering [5-9]. Similar situations hold by studying the phenomena that happen in plasm, in the processes of heat propagation, moisture transfer in soils, the fluid filtration in porous - fissured media and also in the problems of mathematical biology and demography [10].

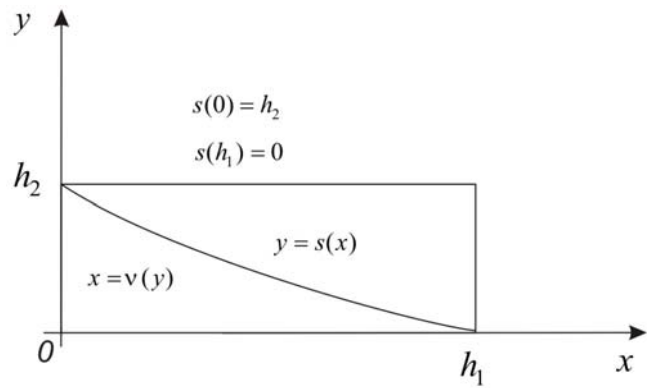

**Fig.1. Graphic representation of the function $s: \overline{G}_1 \times \overline{G}_2$ and its inverse function.**

For equation (1) state the Cauchy problem in the form:
$$(V_{3,j}u)(x) \equiv D_x^3 D_y^j u(x, S(x)) = Z_{3,j}(x), \ x \in G_1,$$
$$(V_{i,2}u)(y) \equiv D_x^i D_y^2 u(v(y), y) = Z_{i,2}(y), \ y \in G_2, \quad (2)$$
$$V_{i,j}u \equiv D_x^i D_y^j u(0, h_2) = Z_{i,j}, \ i = \overline{0,2}; \ j = \overline{0,1},$$

where $S : \overline{G}_1 \to \overline{G}_2$ is an absolutely continuous, strictly decreasing function. $S(0) = h_2$, $S(h_1) = 0$, $v : \overline{G}_2 \to \overline{G}_1$ is a function inverse to the function $y = S(x)$, and the derivative $S_x(x)$ is bounded on $G_1$ (see. fig.1).

Let the conditions $a_{i,j}(x,y) \in L_p(G), i = \overline{0,2}; \ j = \overline{0,1}$ be fulfilled and there exist the functions
$$a_{i,2}^0(x) \in L_p(G_1), i = \overline{0,2},$$
$$a_{3,j}^0(y) \in L_p(G_2), j = \overline{0,1}$$

such that $|a_{i,2}(x,y)| \le a_{i,2}^0(x), |a_{3,j}(x,y)| \le a_{3,j}^0(y)$ almost everywhere on $G$;
$$Z = (Z_{3,2}, Z_{3,0}, Z_{3,1}, Z_{0,2}, Z_{1,2}, Z_{2,2}, Z_{0,0}, Z_{1,0}, Z_{2,0}, Z_{0,1}, Z_{1,1}, Z_{2,1}) \in$$
$$\in E_p^{(3,2)} \equiv L_p(G) \times L_p(G_1) \times L_p(G_1) \times$$
$$\times L_p(G_2) \times L_p(G_2) \times L_p(G_2) \times R \times R \times R \times R \times R \times R,$$

where $R$ is a space of real numbers.

Up to now, in all the papers known in references the Cauchy problem for equation (1) was stated and studied in the classic form:
$$\begin{cases} D_x^i D_y^j u(x, S(x)) = Z^{(i,j)}(x), \ x \in G_1, \ i = \overline{0,2}, \ j = \overline{0,1}, \\ D_x^2 D_y^2 u(v(y), y) = Z^{(4)}(y), \ y \in G_2. \end{cases} \quad (3)$$

Show that conditions (2) and (3) are equivalent. Indeed let $u \in W_p^{(3,2)}(G)$ be a solution of problem (1), (2). Show that it is a solution of problem (1), (3) for some $Z^{(i,j)}(x), i = \overline{0,2}, \ j = \overline{0,1}, Z^{(4)}(y)$. For that consider the identities:
$$u(x, S(x)) = u(0, h_2) + \int_0^x D_x u(\tau, S(\tau))d\tau + \int_{h_2}^{S(x)} D_y u(v(\xi), \xi)d\xi, \quad (4)$$

$$D_x u(x, S(x)) = D_x U(0, h_2) + \int_0^x D_x^2 u(\tau, S(\tau))d\tau + \int_{h_2}^{S(x)} D_x D_y u(v(\xi), \xi)d\xi, \quad (5)$$

$$D_y u(x, S(x)) = D_y u(0, h_2) + \int_0^x D_x D_y u(\tau, S(\tau))d\tau + \int_{h_2}^{S(x)} D_y^2 u(v(\xi), \xi)d\xi, \quad (6)$$

$$D_x D_y u(x, S(x)) = D_x D_y u(0, h_2) + \int_0^x D_x^2 D_y u(\tau, S(\tau))d\tau + \int_{h_2}^{S(x)} D_x D_y^2 u(v(\xi), \xi)d\xi, \quad (7)$$

$$D_x^2 u(x, S(x)) = D_x^2 u(0, h_2) + \int_0^x D_x^3 u(\tau, S(\tau))d\tau + \int_{h_2}^{S(x)} D_x^2 D_y u(v(\xi), \xi)d\xi, \quad (8)$$

$$D_x^2 D_y u(x, S(x)) = D_x^2 D_y u(0, h_2) + \int_0^x D_x^3 D_y u(\tau, S(\tau))d\tau + \int_{h_2}^{S(x)} D_x^2 D_y^2 u(v(\xi), \xi)d\xi, \quad (9)$$

As $u \in W_p^{(2,3)}(G)$ is the solution of problem (1), (2), then from (4)-(9) it follows that

$$D_x^2 D_y u(x, S(x)) = Z_{2,1} + \int_0^x Z_{3,1}(\tau)d\tau + \int_{h_2}^{S(x)} Z_{2,2}(\xi)d\xi = \breve{Z}^{2,1}(x),$$

$$D_x^2 u(x, S(x)) = Z_{2,0} + \int_0^x Z_{3,0}(\tau)d\tau + \int_{h_2}^{S(x)} \breve{Z}^{(2,1)}(v(\xi))d\xi = \breve{Z}^{2,0}(x),$$

$$D_x D_y u(x, S(x)) = Z_{1,1} + \int_0^x \breve{Z}^{(2,1)}(\tau)d\tau + \int_{h_2}^{S(x)} Z_{1,2}(\xi)d\xi = \breve{Z}^{(1,1)}(x),$$

$$D_y u(x, S(x)) = Z_{0,1} + \int_0^x \breve{Z}^{(1,1)}(\tau)d\tau + \int_{h_2}^{S(x)} Z_{0,2}(\xi)d\xi = \breve{Z}^{(0,1)}(x),$$

$$D_x u(x, S(x)) = Z_{1,0} + \int_0^x \breve{Z}^{(2,0)}(\tau)d\tau + \int_{h_2}^{S(x)} \breve{Z}^{(1,1)}(v(\xi))d\xi = \breve{Z}^{(1,0)}(x),$$

$$u(x, S(x)) = Z_{0,0} + \int_0^x \breve{Z}^{(1,0)}(\tau)d\tau + \int_{h_2}^{S(x)} \breve{Z}^{(0,1)}(v(\xi))d\xi = \breve{Z}^{(0,0)}(x).$$

This means that $u(x, y)$ is the solution of problem (1), (3) for

$$Z^{(0,0)}(x) = \breve{Z}^{(0,0)}(x),\ Z^{(1,0)}(x) = \breve{Z}^{(1,0)}(x),\ Z^{(0,1)}(x) = \breve{Z}^{(0,1)}(x),\ Z^{(1,1)}(x) = \breve{Z}^{(1,1)}(x),$$

$$Z^{(2,0)}(x) = \breve{Z}^{(2,0)}(x),\ Z^{(2,1)}(x) = \breve{Z}^{(2,1)}(x),\ Z^{(4)}(y) = Z_{2,2}(y).$$

Vice - versa, if $u \in W_p^{(3,2)}(G)$ is the solution of problem (1), (3), then from the identities (4)-(9) we have

$$Z^{(2,1)}(x) = Z^{(2,1)}(0) + \int_0^x D_x^3 D_y u(\tau, s(\tau)) d\tau + \int_{h_2}^{s(x)} Z^{(4)}(\xi) d\xi,$$

$$Z^{(2,0)}(x) = Z^{(2,0)}(0) + \int_0^x D_x^3 u(\tau, s(\tau)) d\tau + \int_{h_2}^{s(x)} Z^{(2,1)}(v(\xi)) d\xi,$$

$$Z^{(1,1)}(v(y)) = Z^{(1,1)}(0) + \int_0^{v(y)} Z^{(2,1)}(\tau) d\tau + \int_{h_2}^{s(x)} D_x D_y^2 u(v(\xi), \xi) d\xi,$$

$$Z^{(0,1)}(v(y)) = Z^{(0,1)}(0) + \int_0^{v(y)} Z^{(1,1)}(\tau) d\tau + \int_{h_2}^{y} D_y^2 u(v(\xi), \xi) d\xi,$$

$$Z^{(1,0)}(x) = Z^{(1,0)}(0) + \int_0^x Z^{(2,0)}(\tau) d\tau + \int_{h_2}^{s(x)} Z^{(1,1)}(v(\xi)) d\xi,$$

$$Z^{(0,0)}(x) = Z^{(0,0)}(0) + \int_0^x Z^{(1,0)}(\tau) d\tau + \int_{h_2}^{s(x)} Z^{(0,1)}(v(\xi)) d\xi.$$

Hence we find that $u(x, y)$ is the solution of problem (1), (2) for

$$Z_{0,0} = Z^{(0,0)}(0),\ Z_{1,0} = Z^{(1,0)}(0),\ Z_{0,1} = Z^{(0,1)}(0),$$

$$Z_{1,1} = Z^{(1,1)}(0),\ Z_{2,0} = Z^{(2,0)}(0),\ Z_{2,1} = Z^{(2,1)}(0),\ Z_{2,2}(y) = Z^{(4)}(y),$$

$$Z_{3,0}(x) = \frac{d}{dx} F^{(1)}(x),\ F^{(1)}(x) = Z^{(2,0)}(x) - \int_{h_2}^{s(x)} Z^{(2,1)}(v(\xi)) d\xi,$$

$$Z_{3,1}(x) = \frac{d}{dx} F^{(2)}(x),\ F^{(2)}(x) = Z^{(2,1)}(x) - \int_{h_2}^{s(x)} Z^{(4)}(\xi) d\xi,$$

$$Z_{0,2}(y) = \frac{d}{dy} F^{(3)}(y), \ F^{(3)}(y) = Z^{(0,1)}(v(y)) - \int_0^{v(y)} Z^{(1,1)}(\tau) d\tau,$$

$$Z_{1,2}(y) = \frac{d}{dy} F^{(4)}(y), \ F^{(4)}(y) = Z^{(1,1)}(v(y)) - \int_0^{v(y)} Z^{(2,1)}(\tau) d\tau.$$

Thus, conditions (2) are equivalent to the classic form Cauchy conditions (3). Therefore, the classic form Cauchy problems (1), (3) and in non-classical treatment (1), (2) are equivalent in the general case. However, the Cauchy problem in non-classical treatment (1), (3) is more natural by statement than problem (1), (3). In the non-classic treatment no additional conditions of "agreement" type are imposed on the right sides

$$Z_{3,2}(x,y) \in L_p(G), \ Z_{3,j}(x) \in L_p(G_1), \ Z_{i,2}(y) \in L_p(G_2),$$

$$Z_{i,j} \in R, i = \overline{0,2}; \ j = \overline{0,1},$$

this enables to expect that the operator of problem (1), (2) gives homeomorphism. However, in case (1), (3) we can't expect it. The matter is that in this case, for the existence of the solution it is necessary that the conditions

$$F^{(1)}(x) \in W_p^{(1)}(G_1), \ F^{(2)}(x) \in W_p^{(1)}(G_1), F^{(3)}(y) \in W_p^{(1)}(G_2)$$

and $F^{(4)}(y) \in W_p^{(1)}(G_2)$, that may be considered as "agreement" conditions connecting the given functions

$$Z^{(0,1)}(x), \ Z^{(1,1)}(x), Z^{(2,0)}(x), \ Z^{(2,1)}(x), \ Z^{(4)}(y)$$

be fulfilled.

This is an advantage of problem (1), (2) compared to non-classic statement of Cauchy problem (1), (3). Note that some problems in non-classical treatment for pseudoparabolic equations were investigated in the papers of the author's [11-14].